# On the honeycomb conjecture and the Kepler problem

Fu-Gao Song[†‡], Francis Austin[§]

[†] College of Electronic Science & Technology, Shenzhen University, P.R. China
E-mail: songfgao@szu.edu.cn
[‡] Shenzhen Key Laboratory of Micro-nano Photonic Information Technology, P.R. China
[§] Department of Applied Mathematics, The Hong Kong Polytechnic University, Hung Hom, Kowloon, Hong Kong, P.R. China

**Abstract.** This paper views the honeycomb conjecture and the Kepler problem essentially as extreme value problems and solves them by partitioning 2-space and 3-space into building blocks and determining those blocks that have the universal extreme values that one needs. More precisely, we proved two results. First, we proved that the regular hexagons are the only 2-dim blocks that have unit area and the least perimeter (or contain a unit circle and have the least area) that tile the plane. Secondly, we proved that the rhombic dodecahedron and the rhombus-isosceles trapezoidal dodecahedron are the only two 3-dim blocks that contain a unit sphere and have the least volume that can fill 3-space without either overlapping or leaving gaps. Finally, the Kepler conjecture can also be proved to be true by introducing the concept of the minimum 2-dim and 3-dim Kepler building blocks.

## 1. Introduction

The honeycomb conjecture is a 2,000-year-old tiling problem that states that amongst all the conceivable shapes of tiles of equal areas that can be used to cover a plane without gaps or overlaps, the hexagonal tiles are the ones that have the shortest perimeter [1, 2]. Many attempts have been made over the years to prove this result with the latest one being published by T. C. Hales in 2001 using a method that is based on computer computations [3]. Closely related to the honeycomb conjecture is the Kepler problem (also called the 3-dimensional sphere packing problem) that was posed by the famous astronomer J. Kepler in 1611. Kepler conjectured that amongst all the possible sphere packings in $E^3$, the face-centered-cubic (f.c.c.) packing is the one that has the greatest density. While attempting to prove this conjecture in 1694, Newton discussed with Gregory the maximum number of spheres that can touch a certain fixed central sphere. Newton believed that the answer was twelve while Gregory thought that thirteen may be possible. The exact number was finally proved to be twelve by Schütte-van der Waerden in 1953 [4, 5]. Finally, Gauss [6] proved that of all the possible lattice type sphere packings, the f.c.c. packing is indeed the densest one with density $\pi/\sqrt{18}$ in 1831.

In 1883, Barlow [7] discovered the existence of infinitely many other *non-lattice* type packings whose densities are also equal to $\pi/\sqrt{18}$. Barlow constructed these packings by stacking up different layers of planar hexagonal patterns in some pre-determined sequence. Subsequently, Hilbert included the sphere packing problem as the second part of the 18th problem in his famous set of "Hilbert Problems" [8] [9]. In 1947, L. F. Tóth devised one of the first convincing proofs of the two-dimensional analogue of the Kepler problem and predicted in 1965 that computers would one day produce a proof of Kepler's conjecture. In 1993, Wu-Yi Hsiang announced that he had found the best way of proving the Kepler conjecture [10] although it was not believed to be valid by most mathematicians in the field of sphere-packing, e.g., see also [11, 12]. The most recent proof was announced in 1997 by T. C. Hales and it involved the checking of thousands of separate cases by computers. Hales was invited to submit his manuscript to the *Annals of Mathematics* where a committee of 12 experts was appointed to review the paper. Four years later, the committee announced that they were still unable to accept that the proof was correct although no errors were found in the proof.

This paper presents an analytic proof of both the honeycomb conjecture and the Kepler problem by regarding them essentially as extreme value problems.

**2. Some theorems on polygons**

**Definition 1**. A building block is called a *honeycomb building block* if it satisfies the following three conditions: (i) it has unit area, (ii) it has the least perimeter and (iii) these blocks can fill the whole plane without either overlapping or leaving gaps.

**Definition 2**. A building block is called a *2-dim Kepler building block* if it satisfies the following three conditions: (i) it contains a unit circle, (ii) it has the least area and (iii) these blocks can fill the whole plane without either overlapping or leaving gaps.

The honeycomb conjecture and the 2-dim Kepler problem are both extreme value problems. Essentially, we seek a honeycomb building block (respectively, a 2-dim Kepler building block) that satisfies the extreme value conditions and an extra condition. The problems can be solved by elementary extremum methods.

First, we tile the plane by using polygons. It is clear that concave polygons cannot be used for either the honeycomb conjecture or the 2-dim Kepler problem (see Fig. 1, for example, in which $A_1A_2'A_3A_4A_5A_6A_7A_8$ and $A_1A_2A_3A_4A_5A_6A_7A_8$ are two octagons that have the same perimeters but different areas; more precisely, a larger circle may possibly be inscribed inside the convex octagon than inside the concave one). On the other hand, a closed plane curve can be regarded as a polygon having an infinite number of edges. In this sense, therefore, we need only consider convex polygon tilings of 2-space for both the honeycomb conjecture and the 2-dim Kepler conjecture.

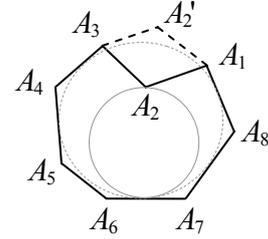

**Fig. 1**

Consider a circle inscribed inside a polygon (see Fig. 2) where $AB$ is not tangent to the circle but $A'B'$ is, and $A'B'//AB$. If we replace the edge $AB$ by $A'B'$, then we obtain a new polygon whose area and perimeter are always less than those of the original one. For the purposes of proving the honeycomb conjecture and the 2-dim Kepler problem, therefore, only circumscribed polygon tilings need be considered.

**Theorem 1.** Suppose that $A_1A_2...A_i...A_n$ is a circumscribed polygon of a circle of radius $r$, where $n \geq 3$ is the number of sides. Among all such polygons of $n$-sides, only the regular polygon has both the least area and the least perimeter. Let $S_n^0(r)$ and $L_n^0(r)$ be the area and the perimeter of the regular polygon, respectively. Then, one has $S_n^0(r) = nr^2 \tan\dfrac{\pi}{n}$, $L_n^0(r) = 2nr\tan\dfrac{\pi}{n}$.

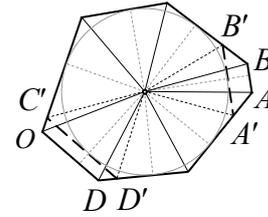

**Fig. 2**

*Proof.* Let $\angle A_iOB = \angle A_iOC = \alpha_i$ at any vertex $A_i$ of the polygon (see Fig. 3). Then the area of the quadrilateral $OBA_iC$ is $r^2\tan\alpha_i$, and $A_iB + A_iC = 2r\tan\alpha_i$. Hence, the area $S_n(r)$ and the perimeter $L_n(r)$ of polygon $A_1A_2...A_i...A_n$ are

$$S_n(r) = r^2 \sum_{1 \le i \le n} \tan \alpha_i,$$

$$L_n(r) = 2r \sum_{1 \le i \le n} \tan \alpha_i = \frac{2S_n(r)}{r}.$$

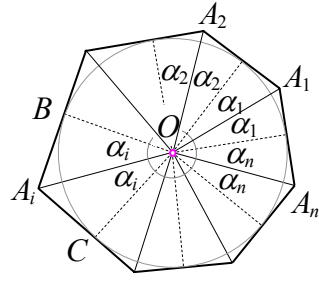

Since $2\sum \alpha_i = 2\pi$, we have $\alpha_n = \pi - (\alpha_1 + \cdots + \alpha_{n-1})$. Hence,

$$S_n(r) = r^2 \sum_{1 \le i \le n-1} \tan \alpha_i - r^2 \tan(\alpha_1 + \cdots + \alpha_{n-1}),$$

$$L_n(r) = 2r \sum_{1 \le i \le n-1} \tan \alpha_i - 2r \tan(\alpha_1 + \cdots + \alpha_{n-1}).$$

**Fig. 3**

The extreme values of $S_n(r)$ and $L_n(r)$ are determined by the equations

$$\frac{\partial S_n(r)}{\partial \alpha_i} = r^2 \left( \frac{1}{\cos^2 \alpha_i} - \frac{1}{\cos^2(\alpha_1 + \cdots + \alpha_{n-1})} \right) = 0, \quad i = 1, \cdots, n-1.$$

And so we must have

$$\cos^2 \alpha_1 = \cos^2 \alpha_2 = \cdots = \cos^2 \alpha_{n-1} = \cos^2(\alpha_1 + \cdots + \alpha_{n-1}) = \cos^2 \alpha_n,$$

i.e.

$$\cos \alpha_i = \pm \cos \alpha_j, \quad i \ne j.$$

Therefore, we must have either $\alpha_i = \alpha_j$ or $\alpha_i = \pi - \alpha_j$, $i \ne j$.

Now, we prove that $\alpha_i \ne \pi - \alpha_j$ for all $i \ne j$. Suppose this is not true. Then we have $\alpha_i + \alpha_j = \pi$, and so it follows from the sum $\sum \alpha_i = \pi$ that $\alpha_k = 0$ for all $k \ne j$ and $k \ne i$. Since this implies that the polygon has only two edges, which is impossible, we must have $\alpha_1 = \alpha_2 = \cdots = \alpha_n$. We have thus shown that the regular $n$-sided polygon is the unique $n$-sided polygon with the property that both $S_n(r)$ and $L_n(r)$ are extreme values. It can also be easily proved that these extreme values are indeed minima, and so we have

$$S_n^0(r) = nr^2 \tan \frac{\pi}{n}, \quad L_n^0(r) = 2nr \tan \frac{\pi}{n}.$$

**Theorem 2.** The following inequalities hold: $S_{n+1}^0(r) \le S_n^0(r)$, $L_{n+1}^0(r) \le L_n^0(r)$, where the equality holds only when $n \to \infty$.

*Proof.* This is because

$$S_n^0(r) = nr^2 \tan\frac{\pi}{n} = \pi r^2 \left(1 + \frac{\pi^2}{3n^2} + \cdots \right),$$

$$S_{n+1}^0(r) = \pi r^2 \left(1 + \frac{\pi^2}{3(n+1)^2} + \cdots \right);$$

$$S_n^0(r) - S_{n+1}^0(r) = \pi r^2 \left(\frac{\pi^2}{3}\frac{2n+1}{n^2(n+1)^2} + \cdots \right) \geq 0,$$

$$L_n^0(r) - L_{n+1}^0(r) = 2\pi r \left(\frac{\pi^2}{3}\frac{2n+1}{n^2(n+1)^2} + \cdots \right) \geq 0.$$

### 3. The hexagonal honeycomb conjecture

The hexagonal honeycomb conjecture states that regular hexagons provide the optimal way to tile the plane into unit areas that have the least perimeter. In general, there are infinitely many ways to tile the plane into unit areas, for instance, by using convex polygons or convex closed curves or any arbitrary closed plane figure. From the proof of Theorem 1, however, we know that all these tilings can be changed to the polygon tilings among which the regular polygon tiling is the optimal. The central point in the honeycomb conjecture is therefore to determine those suitable regular polygons of unit area and least perimeter that will tile the whole plane without leaving gaps or overlapping. The required polygons can have either a finite or an infinite number of edges.

Now, we know that the regular polygon is the one that has the least perimeter amongst all the polygons that have the same number of sides (Theorem 1) and that the perimeters of the regular polygons that enclose a constant area is a decreasing function of its number of sides (Theorem 2). Clearly, every side of the polygons that tile the plane is a common side of two tiling polygons and an interior angle of a regular polygon is $180° \times (n-2)/n$. In order to ensure that the polygons can fill the whole plane without leaving gaps or producing overlaps, therefore, the number of sides $n$ of the regular polygon that can be used must satisfy the equation

$$k\frac{180°(n-2)}{n} = 360°$$

for some integer variable $k$. By Theorem 2, therefore, a proof of the honeycomb conjecture is obtained by determining of the greatest integer $n$ for which the above equation holds. Indeed, the above equation has three solutions

$$k=6, n=3;\ k=4, n=4;\ \text{and}\ k=3, n=6,$$

from which $n=6$ (and $k=3$) is precisely the statement of the honeycomb conjecture. The honeycomb building blocks should therefore be regular hexagons, as claimed by the conjecture.

## 4. The Kepler problem in $E^2$

The Kepler problem in 2-space is also a *space tiling* problem. Specifically, the question is this: which is the polygon with the least area and circumscribing a unit circle that can tile the plane. The answer is exactly the same as that for the hexagonal honeycomb conjecture; the closest packing of unit circles is the regular hexagon packing, where the packing density is equal to

$$\frac{\pi}{S_6^0(1)} = \frac{\pi}{\sqrt{12}}.$$

This conclusion is unique because the 2-dim Kepler building block is unique.

The Kepler conjecture in $E^2$ can also be simply proved by introducing the concept of the *minimum 2-dim Kepler building block*.

**Definition 3.** A building block is called the *minimum 2-dim Kepler building block* if it satisfies the following four conditions: (i) each vertex is located at the center of a unit circle; (ii) the parts of the unit circle cut out by an integral number of blocks precisely form a circle; (iii) it has the least area; (iv) such blocks tile the plane without either overlapping or leaving gaps.

It is clear that the minimum 2-dim Kepler building block is an equilateral triangle of side length 2 (see $\triangle ABC$ in Fig. 4), in which the parts of the unit circle cut out by two blocks precisely form a circle.

The following theorem for the minimum 2-dim Kepler building block is evident.

**Theorem 3.** Let $S_\triangle$ be the area of $\triangle ABC$ in the minimum 2-dim Kepler building block and let $S_c$ be the total area of the parts of the circles cut out by $\triangle ABC$. Define the local packing density for the minimum 2-dim Kepler building block to be $S_c/S_\triangle$. Then $S_c/S_\triangle = \pi/\sqrt{12}$.

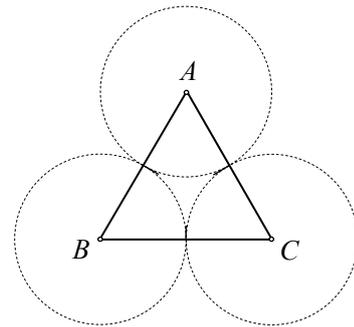

**Fig. 4**

**Definition 4.** The local packing density of a building block is said to be *universal* if the block can be used to fill the whole space.

Obviously, the packing density $\pi/\sqrt{12}$ of the minimum 2-dim Kepler building block is the maximum universal local packing density because the three unit circles in question have already attained the closest packing; no other packings exist that have a packing density (the global packing density included) greater than $\pi/\sqrt{12}$.

In other words, none of the infinitely many packing methods has a packing density greater than $\pi/\sqrt{12}$.

**5. The polygonal pyramids with the maximum local packing density**

**Lemma 1.** Let $O$ be the center of a unit sphere and $OH = 1$ be the vertical of $\triangle ABH$ (see Fig. 5). Let $\theta = \angle AHB$ be constant, $HB \perp AB$, and $BH = x$ be a variable. Let $V_s$ be the partial volume of the sphere cut out by the triangular pyramid $OABH$ and $V_{OABH}$ be that of $OABH$. Then the ratio $V_s / V_{OABH}$ is a decreasing function of the variable $x$.

*Proof.* Clearly, we have $OB \perp AB$, $AB = x\tan\theta$, $AH = x/\cos\theta$ and $OB = (1+x^2)^{1/2}$. Let

$$\varphi = \angle AOB = \arctan[x\tan\theta/(1+x^2)^{1/2}],$$
$$\alpha = \angle BOH = \arctan x,$$
$$\beta = \angle AOH = \arctan(x/\cos\theta);$$
$$p = (\alpha+\beta+\varphi)/4, \; p_1 = (\alpha+\beta-\varphi)/4,$$
$$p_2 = (\alpha-\beta+\varphi)/4, \; p_3 = (-\alpha+\beta+\varphi)/4;$$

then

$$V_s = 4\arctan[(\tan p \tan p_1 \tan p_2 \tan p_3)^{1/2}]/3,$$
$$V_{OABH} = x^2 \tan\theta/6.$$

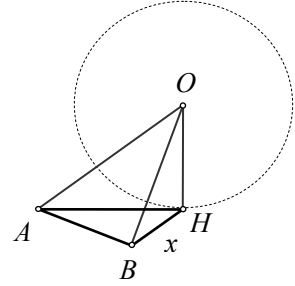

**Fig. 5**

Therefore, $V_s$ and $V_{OABH}$ are both monotonically increasing functions of $x$, and $V_s \to$ constant and $V_{OABH} \to \infty$ as $x \to \infty$. Let $\eta = V_s/V_{OABH}$. Then $\eta$ takes the maximum value $\eta = 1$ when $x = 0$ and tends to the greatest lower bound of 0 (i.e. $\eta \to 0$) as $x \to \infty$. In general, for $0 < x < \infty$ and $0 < \theta < \pi/2$, one has

$$\frac{d\eta}{dx} < 0.$$

Hence, $\eta$ is a decreasing function of $x$.

**Lemma 2.** Let $O$ be the center of a unit sphere and $OH = 1$ be the vertical of $\triangle ABH$ (see Fig. 6). Let $HC \perp AB$, $\angle AHB = \theta$ and $CH = h$ be constants, and $\angle AHC = x$ be a variable. Let $\varpi$ and $V_s$ be respectively the partial solid angle and volume of the sphere cut out by the triangular pyramid $OABH$ and $V_{OABH}$ be the volume of $OABH$. Then $\varpi$, $V_s$ and $V_{OABH}$ take their minimum values and the ratio $\eta = V_s/V_{OABH}$ takes the maximum value when $x = \theta/2$.

*Proof.* It is clear that $\angle BHC = \theta - x$. One also has $AC = h\tan x$, $BC = h\tan(\theta-x)$,

and $AB = h[\tan x + \tan(\theta - x)]$. Hence
$$V_{OABH} = h^2[\tan x + \tan(\theta - x)]/6.$$
On the other hand, we have $OC = (1+h^2)^{1/2}$, $AH = h/\cos x$, $BH = h/\cos(\theta - x)$,
$$\angle AOC = \arctan(AC/OC) = \arctan[h\tan x/(1+h^2)^{1/2}],$$
$$\angle BOC = \arctan[h\tan(\theta-x)/(1+h^2)^{1/2}].$$

Let
$$\alpha = \angle AOH = \arctan(h/\cos x),$$
$$\beta = \angle BOH = \arctan[h/\cos(\theta - x)],$$
$$\varphi = \angle AOB = \angle AOC + \angle COB$$
$$= \arctan\frac{h\tan x}{\sqrt{1+h^2}} + \arctan\frac{h\tan(\theta-x)}{\sqrt{1+h^2}}$$
$$= \arctan\frac{h\sqrt{1+h^2}\sin\theta}{h^2\cos\theta + \cos x\cos(\theta-x)};$$

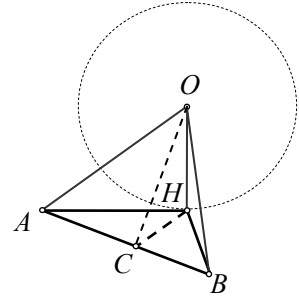

**Fig. 6**

$$p = (\alpha+\beta+\varphi)/4,\ p_1 = (\alpha+\beta-\varphi)/4,$$
$$p_2 = (\alpha-\beta+\varphi)/4,\ p_3 = (-\alpha+\beta+\varphi)/4.$$

Thus
$$\varpi = 4\arctan\left(\sqrt{\tan p\tan p_1\tan p_2\tan p_3}\right),$$
$$V_s = \frac{4}{3}\arctan\left(\sqrt{\tan p\tan p_1\tan p_2\tan p_3}\right),$$
and
$$\eta = \frac{8}{h^2}\frac{\arctan\left(\sqrt{\tan p\tan p_1\tan p_2\tan p_3}\right)}{\tan x + \tan(\theta-x)}.$$

Finally, since $\varpi$, $V_s$, $V_{OABH}$ and $\eta$ are all symmetric functions about $x = \theta/2$, one has
$$\frac{dV_{OABH}}{dx} = 0,\ \frac{dV_s}{dx} = 0,\ \frac{d\varpi}{dx} = 0\ \text{and}\ \frac{d\eta}{dx} = 0$$
when $x = \theta/2$. It can be easily verified that $\varpi$, $V_s$ and $V_{OABH}$ take the minimum values and $\eta$ the maximum at $x = \theta/2$, and so the result follows.

**Definition 5.** A polygonal pyramid is called a *tight polygonal pyramid* if it satisfies the following conditions: (i) its vertex is at the center of a unit sphere; (ii) if another unit sphere is put inside the polygonal pyramid, then it is tangent to the unit sphere and all the lateral surfaces of the polygonal pyramid.

The following lemma is often used for the tight polygonal pyramid.

**Lemma 3.** Let the base surface $A_1A_2\ldots A_n$ of a tight polygonal pyramid $OA_1A_2\ldots A_n$ be tangent to the unit sphere $O$ at point $H$. Then the distances from $H$ to each side of

$A_1A_2…A_n$ are all equal to $1/\sqrt{3}$.

*Proof.* In Fig. 7, $HH_1 \perp A_1A_2$ and $OD_1$ is the extension line of $OH_1$ that is tangent to the unit sphere $O'$ at point $D_1$. Since $OH = 1$, $O'D_1 = 1$, $OO' = 2$, one has $\angle D_1OO' = 30°$ and hence $HH_1 = 1/\sqrt{3}$.

The imaginary circle of radius $1/\sqrt{3}$ on the tangent plane $A_1A_2…A_n$ in Fig. 7 is a section of the cone formed from the sphere $O'$ with vertex $O$. The same imaginary circles appear in Figures 8, 10, 11 and 12 (*a*) and (*b*) below.

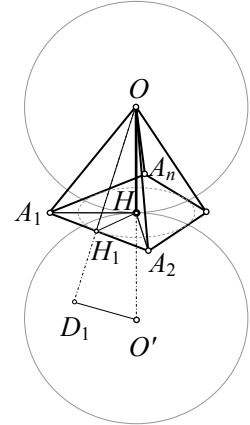

**Fig. 7**

### 6. The thirteen spheres problem

Here, we consider incidentally the thirteen spheres problem.

Let $\varpi_n$ and $\tau_n$ denote respectively the solid angle and the volume of the unit sphere that is cut out by a regular tight polygonal pyramid of *n* edges and unit altitude. Let $V_n$ be its volume and $\eta_n = \tau_n / V_n$. We have the following theorem.

**Theorem 4.** Of all the tight polygonal pyramids of *n* edges, only the regular tight polygonal pyramids have the minimum $\varpi_n$, $\tau_n$ and $V_n$ and the maximum $\eta_n$.

*Proof.* In the triangular pyramid $OA_1HA_2$ (see Fig. 8), $\varpi$, $V_s$ and $V_{OABH}$ all take their minimum values and $\eta$ takes the maximum when $\angle A_1HH_1 = \angle A_2HH_1$ (Lemma 1). In the tight polygonal pyramid $OA_1A_2…A_n$, $h = HH_1 = HH_2 = … = HH_n = 1/\sqrt{3}$ (Lemma 2). Hence, $\Delta A_1HH_1 \equiv \Delta A_2HH_1 \equiv \Delta A_2HH_2 \equiv … \equiv \Delta A_1HH_n$, and so the tight polygonal pyramid $OA_1A_2…A_n$ is in fact a regular tight polygonal pyramid and the result follows.

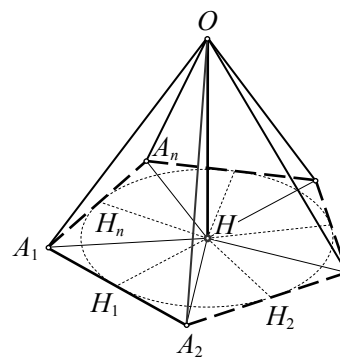

**Fig. 8**

The hexagonal pyramids are of no use in the 3-dim Kepler problem because their bases they do not form a polyhedron. Here, we only study the pentagonal pyramids, the quadrilateral pyramids and the triangular pyramids.

Table 1 shows the values of $\varpi_n$, $\tau_n$, $V_n$, $\eta_n$ and $4\pi/\varpi_n$ for the regular tight pentagonal pyramid ($n = 5$), the regular tight quadrilateral pyramid ($n = 4$) and the regular tight triangular pyramid ($n = 3$).

**Table 1. The values of $\varpi_n$, $\tau_n$, $V_n$, $\eta_n$ and $4\pi/\varpi_n$ for different $n$**

| $n$ | 3 | 4 | 5 |
|---|---|---|---|
| $\varpi_n$ | 1.194812833 | 1.010721021 | 0.9425295571 |
| $\tau_n$ | 0.3982709444 | 0.3369070068 | 0.3141765190 |
| $V_n$ | $\sqrt{3}/3$ | $4/9$ | $5\tan 36°/9$ |
| $\eta_n$ | 0.6898255109 | 0.7580407654 | 0.7783683853 |
| $4\pi/\varpi_n$ | 10.51743860 | 12.43307536 | 13.33260110 |

The thirteen spheres problem has been solved by K. Schütte and B. L. van der Waerden in 1953 [4]. Here, we show in a simpler way why thirteen unit spheres cannot simultaneously touch one fixed unit sphere.

**Definition 6.** A polygonal pyramid is called a *quasi-tight polygonal pyramid* if it satisfies the following conditions: (i) its vertex is at the center of a unit sphere; (ii) if another unit sphere is put inside the polygonal pyramid, then it is tangent to the unit sphere and only some lateral surfaces of the polygonal pyramid and have no intersection with the other lateral surfaces.

The following lemma is a direct consequence of Lemma 2.

**Lemma 4.** Let $\omega_n$ ($n = 3, 4, 5, \ldots$) be the partial solid angles of the sphere cut out by a tight polygonal pyramid. Then $\omega_n \geq \varpi_n$ ($n = 3, 4, 5, \ldots$).

The following lemma is also a direct consequence of Lemma 1 and Lemma 2.

**Lemma 5.** Let $\omega'_n$ ($n = 3, 4, 5, \ldots$) be the partial solid angles of the sphere cut out by a quasi-tight polygonal pyramid. Then $\omega'_n \geq \omega_n \geq \varpi_n$ ($n = 3, 4, 5, \ldots$).

It is clear from Table 1 that a unit sphere cannot be partitioned by 13 tight quadrilateral- (or triangular-) pyramids without overlapping because $13\omega_n \geq 13\varpi_n > 4\pi$ ($n = 3, 4$). A unit sphere, however, can be partitioned by at most 12 tight or quasi-tight quadrilateral pyramids, or by at most 10 tight or quasi-tight triangular pyramids. Similarly, a unit sphere cannot be partitioned by 13 pentagonal pyramids even though $13\varpi_5 < 4\pi$ because a tridecahedron cannot be built from thirteen pentagons (a polyhedron can only possibly be built from an even number of pentagons, e.g., a dodecahedron).

Furthermore, even though a tridecahedron can be built from eleven pentagons together with a quadrilateral and a triangle (see Fig. 9(*a*)), a unit sphere cannot be partitioned without overlapping by using eleven tight pentagonal pyramids together with a tight quadrilateral pyramid and a tight triangular pyramid because we have $11\omega_5 + \omega_4 + \omega_3 \geq 11\varpi_5 + \varpi_4 + \varpi_3 > 4\pi$.

Again, even though a tridecahedron can be built from ten pentagons and three quadrilaterals as shown in Fig. 9(*b*), a unit sphere cannot be partitioned without overlapping by ten quasi-tight pentagonal pyramids and three quasi-tight quadrilateral pyramids although $10\varpi_5 + 3\varpi_4 < 4\pi$. This is due to the many gaps that must exist in every arrangement of the ten pentagonal pyramids. Take for example the case of five pentagonal pyramids surrounding a single pentagonal pyramid. In this case, at least two of the pentagonal pyramids must be quasi-tight. Likewise, the three quadrilateral pyramids are also not the regular tight quadrilateral pyramids but the "quasi-tight" ones. Therefore, we must have $\Sigma\omega'_5 + \Sigma\omega'_4 > 4\pi$.

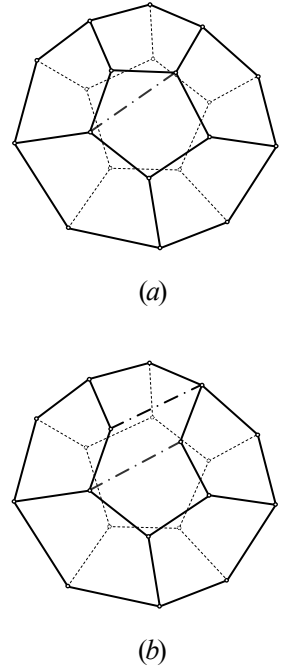

(*a*)

(*b*)

Fig. 9

In short, it is impossible for thirteen unit spheres to touch a single unit sphere simultaneously.

## 7. The Kepler conjecture in $E^3$

**Definition 7.** A building block is called a *3-dim Kepler building block* if it satisfies the following three conditions: (i) it contains a unit sphere, (ii) it has the least volume and (iii) such blocks can fill the whole 3-space without either overlapping or leaving gaps.

**Definition 8**. A sphere packing is called the *tight sphere packing* if every unit sphere has twelve circumscribed unit spheres. These are also called the *periphery spheres* of the *central sphere*.

The following lemma is evident.

**Lemma 6.** In the tight sphere packing, every periphery sphere is tangent to, except for the central sphere, those spheres that belong to its own periphery spheres existing in the other eleven periphery spheres.

**Remark.** In fact, we can see below that each one of the 12 periphery spheres of a given central sphere has in turn its own set of 12 periphery spheres that is made up of the following: the original "mother" central sphere, 4 periphery spheres of the "mother" central sphere (see Fig. 12 (c) and (d) below), and 7 "outer layer" spheres.

Since a unit sphere has at most twelve circumscribed unit spheres, the tight sphere packing is the closest packing. The 3-dim Kepler problem is therefore equivalent to determining the upper bound of the packing densities in the tight sphere packing. Although not the simplest method, the method of tight polygonal pyramids is still effective in solving the 3-dim Kepler problem.

**Lemma 7.** Each sphere in the tight sphere packing has its own definite circumscribed dodecahedron.

*Proof.* Since every unit sphere has twelve periphery spheres in the tight sphere packing, the tangent planes to the central sphere that touch the periphery spheres form a circumscribed dodecahedron.

The circumscribed dodecahedron above will be called the *tight dodecahedron*.

It is clear from Definition 8 and Lemma 7 that the 3-dim Kepler building blocks are tight dodecahedrons.

**Definition 9.** Let the centers of every pair of mutually tangent periphery spheres be joined together in the tight sphere packing to obtain a polyhedron. Then the collection of those tangent planes on the tight dodecahedron whose tangent points are determined by the vertices of a definite face of the above polyhedron is called a *tight ring* of the tight dodecahedron.

**Lemma 8.** The tangent planes of a tight ring of the tight dodecahedron have a common vertex and the plane angles at the common vertex are all identical.

*Proof.* The conclusion is evident if the number of tangent planes is three. We prove the conclusion to be true for four tangent planes below.

In Fig. 10, the four tangent planes of a tight ring are expressed respectively by four tangent points $O_1$, $O_2$, $O_3$ and $O_4$. Since the two periphery spheres whose tangent planes have a common side are tangent to each other, and $O_1A \perp A_1A$, $O_4A \perp A_1A$, $O_1B \perp B_1B$, $O_2B \perp B_1B$, $O_2C \perp C_1C$, $O_3C \perp C_1C$, $O_3D \perp D_1D$, $O_4D \perp D_1D$, then we have $O_1A = O_4A = O_1B = O_2B = O_2C = O_3C = O_3D = O_4D = 1/\sqrt{3}$ (Lemma 3), and so the eight right triangles are all identical. Since the four right triangles on the tangent planes $O_1$ and $O_2$ are all identical, the extension lines of $A_1A$, $B_1B$ and $C_1C$ must intersect at a point $V$; likewise, those of $B_1B$, $C_1C$ and $D_1D$ must also intersect at the point $V$. Therefore, the four tangent planes have a common vertex $V$, thus $\angle AVB \equiv \angle BVC \equiv \angle CVD \equiv \angle DVA$ and $\angle AO_1V \equiv \angle BO_1V \equiv \angle BO_2V \equiv \angle CO_2V \equiv \angle CO_3V \equiv \angle DO_3V \equiv \angle DO_4V \equiv \angle AO_4V$.

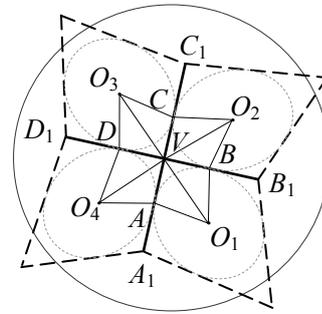

Fig. 10

**Remark.** Although three distinct tangent planes always intersect in a unique common vertex, four distinct tangent planes however may intersect at two vertices. Consider for example Fig. 10, where the tangent planes $O_1$, $O_2$ and $O_3$ intersect at $V_1$ and tangent planes $O_1$, $O_3$ and $O_4$ intersect at $V_2$, which may be distinct from $V_1$. Lemma 8, however, guarantees that $V_1$ and $V_2$ are the same for four tangent planes of a tight ring of the tight dodecahedron. This also implies that every face of a tight dodecahedron is the base of the tight polygonal pyramid.

**Lemma 9.** In the tight sphere packing, every sphere can be partitioned by twelve tight polygonal pyramids without overlapping or leaving gaps.

*Proof.* This partition can be attained as long as we make a plane through every edge of the tight dodecahedron and the center of the central sphere.

Therefore, the last step in solving the Kepler problem is to seek the tight polygonal pyramids with the minimum base area such that the tight dodecahedron built from the bases has the minimum volume.

**Definition 10.** The packing density of a 3-dim building block is defined to be the ratio of the sum of the partial volumes of the unit spheres cut out by the block to its own volume.

It was shown in Table 1 that it is not possible to partition a unit sphere by regular tight polygonal pyramids because the ratio $4\pi / \varpi_n$ is not an integer. A dodecahedron, on the other hand, can be built from twelve pentagons or quadrilaterals. Indeed, the regular dodecahedron consists of twelve regular pentagons, even though the corresponding regular pentagonal pyramids are not the "tight" ones. The packing density of the regular dodecahedron is therefore of no use because it is unsuitable for other unit spheres. If we partition the unit sphere by using twelve tight polygonal (or pentagonal, or quadrilateral) pyramids, then it might be possible to obtain a packing density which is suitable for all unit spheres. We now try this partition method for which the following lemma is significant.

**Lemma 10.** If the bases of two tight polygonal pyramids have the same composition of interior angles, then their base areas and the respective partial solid angles cut out from the unit sphere are all the same.

*Proof.* We prove the conclusion only for the tight quadrilateral pyramid. Let $OA_1A_2A_3A_4$ be a tight quadrilateral pyramid (See Fig. 11), where $O$ is the center of a unit sphere, $OH = 1$ is the vertical of $A_1A_2A_3A_4$, $HH_1 \perp A_1A_2$, $HH_2 \perp A_2A_3$, $HH_3 \perp A_3A_4$ and $HH_4 \perp A_4A_1$, and $HH_1 = HH_2 = HH_3 = HH_4 = 1/\sqrt{3}$. The area $S$ of the quadrilateral is then equal to $S = (\tan\alpha_1 + \tan\alpha_2 + \tan\alpha_3 + \tan\alpha_4)/3$, which depends only on the values of $\alpha_1$, $\alpha_2$, $\alpha_3$ and $\alpha_4$. Likewise, the solid angle cut out from the unit sphere by the tight quadrilateral pyramid $OA_1A_2A_3A_4$ is equal to the sum of the solid angles cut out by the quadrilateral pyramids $OA_1H_1HH_4$, $OA_2H_2HH_1$, $OA_3H_3HH_2$ and $OA_4H_4HH_3$, which also depends only on the values of $\alpha_1$, $\alpha_2$, $\alpha_3$ and $\alpha_4$. Therefore, the lemma holds.

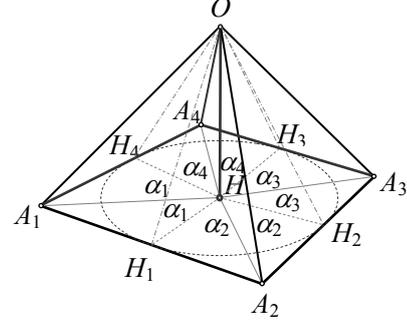

Fig. 11

We now seek the tight polygonal pyramids that have the minimum base area.

Since a tight quadrilateral pyramid (e.g. $OA_1A_2A_3A_4$ in Fig. 11) is characterized by four variables $\alpha_1$, $\alpha_2$, $\alpha_3$ and $\alpha_4$ and one constraint $\alpha_1 + \alpha_2 + \alpha_3 + \alpha_4 = \pi$, it can be completely described by only three independent variables. For convenience, let us *assume*, without loss of generality (see below), that the compositions of the interior angles of the bases of the twelve tight polygonal pyramids are the same. It follows from Lemma 10 that these twelve tight polygonal pyramids have the same volume, solid angles and base area. As the twelve bases form a dodecahedron, every solid angle is $\pi/3$. Also, since the base areas of the twelve tight polygonal pyramids are dependent entirely on $\alpha_1$, $\alpha_2$, $\alpha_3$ and $\alpha_4$, we can select these values such that the base areas are minimum and the packing density of the dodecahedron maximum. This is clearly an extreme value problem.

The area of the quadrilateral $A_1A_2A_3A_4$ in Fig. 11 is given by
$$s = (\tan\alpha_1 + \tan\alpha_2 + \tan\alpha_3 + \tan\alpha_4)/3, \quad \alpha_4 = \pi - \alpha_1 - \alpha_2 - \alpha_3,$$
where $\alpha_1$, $\alpha_2$ and $\alpha_3$ are regarded as the independent variables. These values must also be chosen such that the solid angle at the vertex $O$ is $\pi/3$. There are hence only two independent variables for this problem.

The extreme value of $s$ can be determined in at least three ways. First, let
$$\frac{\partial s}{\partial \alpha_1} = 0, \quad \frac{\partial s}{\partial \alpha_2} = 0, \quad \frac{\partial s}{\partial \alpha_3} = 0.$$

One then obtains $\alpha_1 = \alpha_2 = \alpha_3 = \alpha_4 = \pi/4$ and $OA_1A_2A_3A_4$ is a regular tight quadrilateral pyramid in this case. As discussed earlier, a dodecahedron cannot be built from the bases of these quadrilateral pyramids.

A second method is to let

$$\frac{\partial s}{\partial \alpha_1} = 0, \quad \frac{\partial s}{\partial \alpha_2} = 0.$$

One then obtains $\alpha_1 = \alpha_2 = \alpha_4 = \alpha$ and $\alpha_3 = \pi - 3\alpha$, where $\alpha$ is a variable that can be chosen to ensure that the solid angle at the vertex $O$ is $\pi/3$. To achieve this, one must have

$$\alpha_1 = \alpha_2 = \alpha_4 = 51.178151° \text{ and } \alpha_3 = 26.465547°.$$

No dodecahedrons circumscribing a unit sphere can be built from the bases of these tight quadrilateral pyramids because

$$\angle A_2 A_3 A_4 = 2 \times (90° - \alpha_3) = 127.06891°,$$

and every vertex of a polyhedron is made up of at least three plane angles, and

$$3 \angle A_2 A_3 A_4 = 381.20672° > 360°.$$

This implies that the dodecahedron cannot be built from these bases.

The third way is to fix $\alpha_2$ and $\alpha_3$ and let $\alpha_1$ be variable and put

$$\frac{\partial s}{\partial \alpha_1} = 0.$$

One then obtains $\alpha_1 = \alpha_4 = \alpha$ and $s$ can be rewritten as

$$s = (2\tan\alpha + \tan\alpha_2 + \tan\alpha_3)/3, \quad \alpha_3 = \pi - 2\alpha - \alpha_2.$$

Next we fix $\alpha$ and let $\alpha_2$ be variable, and let

$$\frac{\partial s}{\partial \alpha_2} = 0.$$

One then obtains $\alpha_1 = \alpha_4 = \alpha$, $\alpha_2 = \alpha_3 = \pi/2 - \alpha$. In order to ensure that the solid angle at the vertex $O$ is $\pi/3$, we must have either $\alpha = \arctan 2^{1/2}$ or $\alpha = \pi/2 - \arctan 2^{1/2}$. These bases can now precisely form the circumscribed dodecahedron of a unit sphere, and the tight dodecahedron has packing density $\pi/\sqrt{18}$.

Since the area $s$ is independent of the order of $\alpha_1$, $\alpha_2$, $\alpha_3$ and $\alpha_4$, there are in fact two kinds of tight quadrilateral pyramids: one has rhombic bases and the other one isosceles trapezium bases, as shown in Fig. 12 (*a*) and (*b*), in which $\alpha_1 = \arctan 2^{1/2}$ and $\alpha_2 = \pi/2 - \arctan 2^{1/2}$.

Note that twelve rhombuses can make up a dodecahedron, as shown in Fig. 12 (*c*), which corresponds to the well-known f.c.c packing, but twelve isosceles trapeziums cannot. Alternatively, a dodecahedron can also be made up of six isosceles trapeziums

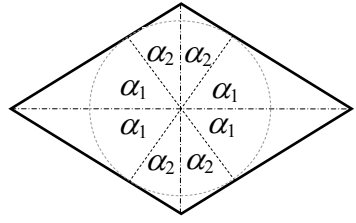
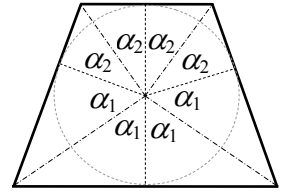

(*a*) Rhombus  (*b*) Isosceles trapezium

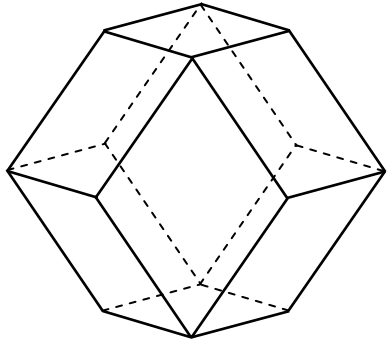
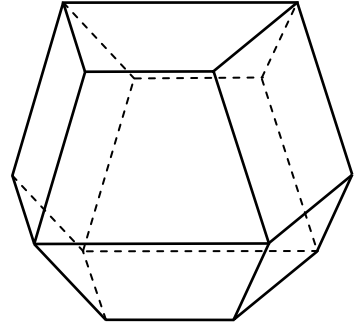

(*c*) Rhombic dodecahedron  (*d*) Rhombus-isosceles trapezoidal dodecahedron and six rhombuses, as shown in Fig. 12 (*d*), which corresponds to the well-known hexagonal close packing. These two kinds of dodecahedrons can be used either alone or together to fill the whole 3-space without either overlapping or leaving gaps. They are, therefore, all the 3-dim Kepler building blocks with packing density $\pi/\sqrt{18}$.

Fig. 12

Now the base area of a tight pentagonal pyramid is given by

$$s = (\tan\alpha_1 + \tan\alpha_2 + \tan\alpha_3 + \tan\alpha_4 + \tan\alpha_5)/3, \quad \alpha_5 = \pi - \alpha_1 - \alpha_2 - \alpha_3 - \alpha_4.$$

Repeating the above discussions, one can prove that no tight pentagonal pyramids exist whose bases can form a circumscribed dodecahedron of a unit sphere.

The following theorem hence follows.

**Theorem 5.** There are two kinds of 3-dim Kepler building blocks with packing density $\pi/\sqrt{18}$; they are the rhombic dodecahedrons and the rhombus-isosceles trapezoidal dodecahedrons.

It should be pointed out that the f.c.c. packing, the hexagonal close packing and the other sphere packings discovered by Barlow [7] can be made up of these two kinds of 3-dim Kepler building blocks.

Nevertheless, Theorem 5 is obtained under the hypothesis that *the compositions of the interior angles of the bases of the twelve tight polygonal pyramids are the same* (disregarding the order of the interior angles). Does this hypothesis hold in general? The answer is positive. We know that every quadrilateral has three independent interior angle variables, and a quadrilateral that has the minimum area has only one which ensures that the solid angle of the corresponding tight quadrilateral pyramid has the suitable value mentioned above. Hence, each one of those twelve quadrilaterals that make up a tight dodecahedron with the minimum volume has only one independent interior angle variable. As Lemma 8 guarantees that all plane angles with a common vertex are identical in a tight dodecahedron, we deduce that all quadrilaterals with a common vertex have the same composition of interior angles and that all the twelve quadrilaterals of the tight dodecahedron have identical interior angles. Therefore, the hypothesis used for deducing Theorem 5 indeed holds and thus Theorem 5 holds in general.

Is there a sphere packing in which the global packing density is greater than $\pi/\sqrt{18}$? The answer is negative. As mentioned above, the tight sphere packing is the closest packing in which every unit sphere is fully partitioned by twelve tight polygonal pyramids (Lemma 9). Theorem 5 has already summed up all the results of the closest packing and there are no others.

Clearly, the Kepler conjecture is true and the tight sphere packing contains the f.c.c. packing, the hexagonal close packing and all the other sphere packings that have been discovered by Barlow. Apart from these, there are no other tight sphere packings.

In fact, there is a simpler way to solve the Kepler problem by introducing the concept of the *minimum 3-dim Kepler building block*.

**Definition 11.** A 3-dim building block is called the *minimum 3-dim Kepler building block* if it satisfies the following four conditions: (i) its vertices are all located at the center of a unit sphere, (ii) it has the minimum volume, (iii) the parts of the unit spheres cut out by an integral number of blocks can precisely make up a whole unit sphere and (iv) such blocks can fill the whole 3-space without either overlapping or leaving gaps.

Since no global packing density can exceed the maximum universal local packing density, the Kepler problem would be completely solved as soon as the minimum 3-dim Kepler building blocks is found and its local packing density determined.

We shall find the minimum 3-dim Kepler building block by verifying each building block, starting from the simplest 3-dim block.

The simplest 3-dim building block is a tetrahedron that is made up of four triangles. However, because tetrahedrons do not satisfy condition (iv) of Definition 11, they are not minimum 3-dim Kepler building blocks.

A pentahedron can be formed by either four triangles and a quadrilateral or two triangles and three quadrilaterals. The former construction is not a minimum 3-dim Kepler building block because it does not satisfy condition (iv) of Definition 11 but the latter may possibly be if it satisfies certain conditions.

Instead of the pentahedron, which is made up of two triangles and three quadrilaterals, let us first discuss the hexahedron which is built from six quadrilaterals that has the closest local sphere packing and is able to fill the whole 3-space without either overlapping or leaving gaps. The conclusion of the former can be readily applied to the latter.

The building blocks that have the closest local sphere packing must have the minimum volume. This fact will be particularly useful when seeking the minimum 3-dim Kepler building blocks.

Consider the tetrahedron with minimum volume and edge lengths not less than 2. This is an extreme value problem. The four vertices of a tetrahedron can be positioned as shown in Fig. 13, in which one adds the following constraints:

$$x \geq 1, y \geq 0, y_3 \geq 0, z \geq 0,$$
$$(x_2-x)^2+y^2 \geq 4, (x_2+x)^2+y^2 \geq 4,$$
$$(x_3-x)^2+y_3^2+z^2 \geq 4, (x_3+x)^2+y_3^2+z^2 \geq 4,$$
$$(x_3-x_2)^2+(y_3-y)^2+z^2 \geq 4.$$

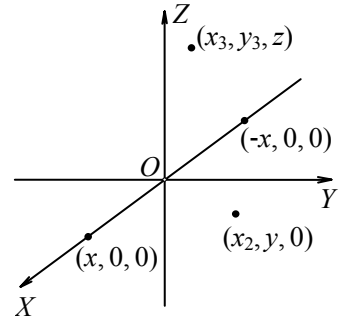

Fig. 13

Note that the constraint $y_3 \geq 0$ is based on the consideration that the tetrahedron would be used for seeking the pentahedrons and the hexahedrons which have the minimum volume. Clearly, the volume of the tetrahedron is

$$V_4 = xyz/3,$$

which is an increasing function of the edge lengths.
One can change the constraints as follows:

$$x \geq 1, y \geq 0, y_3 \geq 0, z \geq 0,$$
$$(x_2-x)^2+y^2 = 4, (x_2+x)^2+y^2 = 4,$$
$$(x_3-x)^2+y_3^2+z^2 = 4, (x_3+x)^2+y_3^2+z^2 = 4,$$
$$(x_3-x_2)^2+(y_3-y)^2+z^2 = 4.$$

Then one obtains $x_2 = x_3 = 0$, $y_3 = (2-x^2)/(4-x^2)^{1/2} \geq 0$ (hence $1 \leq x \leq 2^{1/2}$), $y = (4-x^2)^{1/2}$, and $z = 2[(3-x^2)/(4-x^2)]^{1/2}$ and thus we have

$$V_4 = 2x(3-x^2)^{1/2}/3, \ (1 \leq x \leq 2^{1/2})$$

which attains a maximum value at $x = 1.5^{1/2}$. Similarly, the local minima of $V_4$ are attained at $x = 0$ and $x = 2^{1/2}$ and it turns out that the two minima are identical and

$$V_4 = 2\sqrt{2}/3.$$

There are two kinds of tetrahedron, with vertices located at

$$(-1, 0, 0), (1, 0, 0), (0, 3^{1/2}, 0), (0, 3^{-1/2}, 2/1.5^{1/2})$$

and

$$(-2^{1/2}, 0, 0), (2^{1/2}, 0, 0), (0, 2^{1/2}, 0), (0, 0, 2^{1/2}).$$

These tetrahedrons, however, cannot fill the whole 3-space.

If we now add an identical triangle to the base of the above tetrahedron, then we obtain a pentahedron consisting of four triangles and a quadrilateral whose volume is

$$V_5 = 2xyz/3,$$

with a minimum value of $4\sqrt{2}/3$. There are two kinds of pentahedrons, the five vertices of which are located at

$$(-1, 0, 0), (1, 0, 0), (2, 3^{1/2}, 0), (0, 3^{1/2}, 0) \text{ and } (0, 3^{-1/2}, 2/1.5^{1/2})$$

and

$$(-2^{1/2}, 0, 0), (0, -2^{1/2}, 0), (2^{1/2}, 0, 0), (0, 2^{1/2}, 0) \text{ and } (0, 0, 2^{1/2}).$$

These pentahedrons, however, cannot fill the whole 3-space.

On the other hand, if we let the base and the altitude of the above pentahedrons be those of a hexahedron, then the volume of this hexahedron is

$$V_6 = 2xyz,$$

with a minimum value of $4\sqrt{2}$. There are thus two kinds of hexahedrons, the eight vertices of which are located at

$$(-1, 0, 0), (1, 0, 0), (2, 3^{1/2}, 0), (0, 3^{1/2}, 0); (0, 3^{-1/2}, 2/1.5^{1/2}),$$
$$(2, 3^{-1/2}, 2/1.5^{1/2}), (3, 3^{1/2}+3^{-1/2}, 2/1.5^{1/2}), (1, 3^{1/2}+3^{-1/2}, 2/1.5^{1/2})$$

and

$$(-2^{1/2}, 0, 0), (0, -2^{1/2}, 0), (2^{1/2}, 0, 0), (0, 2^{1/2}, 0); (0, 0, 2^{1/2}),$$
$$(2^{1/2}, -2^{1/2}, 2^{1/2}), (8^{1/2}, 0, 2^{1/2}), (2^{1/2}, 2^{1/2}, 2^{1/2})$$

with their respectively projections in the X-Y plane as shown in Fig. 14 (*a*) and (*b*). It therefore follows that all the edge lengths of these two hexahedrons are 2 and that the eight spheres in question have attained the closest stack.

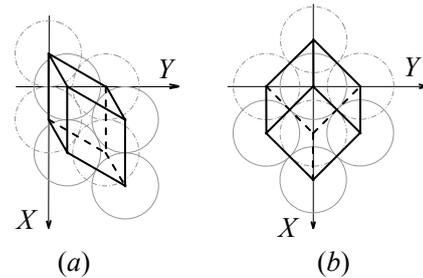

Fig. 14

So there are two different kinds of parallelepipeds, one consisting of six rhombuses and the other consisting of four rhombuses and two squares. In addition, the parts of the unit spheres cut out by the eight vertices do precisely make up one unit sphere with packing density $\pi/\sqrt{18}$. The whole of 3-space can be filled by either one of these two kinds of parallelepipeds or in combination together. Using these two kinds of blocks, we can stack the f.c.c. packing, the hexagonal close packing and the other Barlow sphere packings. These are thus the only tight sphere packings that exist..

Both of the above parallelepipeds can be partitioned along the two parallel diagonals of any two parallel planes into two identical pentahedrons that are made up of two triangles and three quadrilaterals. The pentahedrons are in fact triangular prisms that are all minimum 3-dim Kepler building blocks with a packing density of $\pi/\sqrt{18}$. This is the maximum universal local packing density.

To sum up, we have the following theorems.

**Theorem 6 (Kepler theorem).** No sphere packing in 3-space can have a global packing density that exceeds $\pi/\sqrt{18}$.

**Theorem 7.** The f.c.c. packing, the hexagonal close packing and the sphere packings discovered by Barlow are the only tight sphere packings that exist.